\documentclass[11pt,a4paper,reqno]{amsart}
\usepackage{amsmath,amssymb,amsthm,enumitem,graphicx,color,tikz}
\usetikzlibrary{calc}
\usetikzlibrary{decorations.markings}

\title{Even-integer continued fractions and the Farey tree}
\author{Ian Short}
\author{Mairi Walker}
\address{Department of Mathematics and Statistics\\The Open University\\Milton Keynes MK7 6AA\\ United Kingdom}
\email{ian.short@open.ac.uk}
\email{mairi.walker@open.ac.uk}
\date{\today}

\newtheorem{theorem}{Theorem}[section]
\newtheorem*{theorem*}{Theorem}
\newtheorem{lemma}[theorem]{Lemma}
\newtheorem{corollary}[theorem]{Corollary}
\theoremstyle{definition}

\theoremstyle{plain}

\numberwithin{equation}{section}
\renewcommand{\leq}{\leqslant}
\renewcommand{\geq}{\geqslant}

\definecolor{col2}{RGB}{95,141,211}
\definecolor{col1}{RGB}{215,227,224}
\definecolor{col3}{RGB}{255,255,255}

\newcommand{\horo}[3]{%
   \draw[fill=#3] (#1,#2) circle (#2);
}

\tikzset{middlearrow/.style={
        decoration={markings,
            mark= at position 0.5 with {\arrow{#1}} ,
        },
        postaction={decorate}
    }
}

\begin{document}

\begin{abstract}
Singerman introduced to the theory of maps on surfaces an object that is a universal cover for any map. This object is a tessellation of the hyperbolic plane together with a certain subset of the ideal boundary. The 1-skeleton of this tessellation comprises the edges of an infinite tree whose vertices belong to the ideal boundary. Here we show how this tree can be used to give a beautiful geometric representation of even-integer continued fractions. We use this representation to prove some of the fundamental theorems on even-integer continued fractions that are already known, and we also prove some new theorems with this technique, which have familiar counterparts in the theory of regular continued fractions. 
\end{abstract}

\maketitle

\section{Introduction}

\begin{figure}[ht]
\centering
\includegraphics{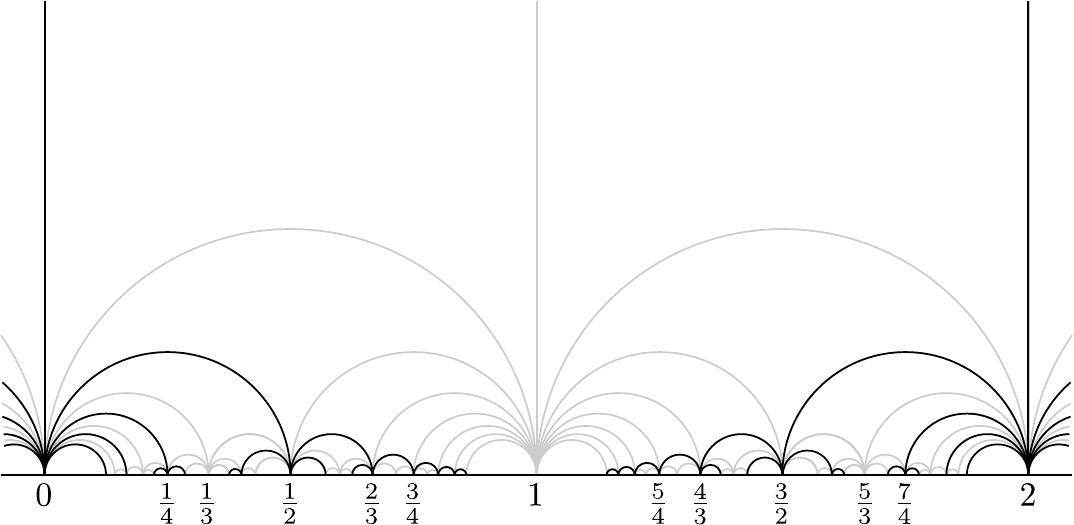}
\caption{The Farey tree superimposed with the Farey graph}
\label{ixe}
\end{figure}

In \cite{Si1988}, Singerman introduced a tessellation of the hyperbolic plane that can be used as a universal cover for any map on a surface. To describe this universal tessellation, we first define the well known \emph{Farey graph}, written as~$\mathcal{G}$. We use the upper half-plane model of the hyperbolic plane, denoted by~$\mathbb{H}$, along with the ideal boundary of~$\mathbb{H}$, which is the extended real line~$\mathbb{R}_\infty$ (that is, the real line~$\mathbb{R}$ with the point~$\infty$ attached). The Farey graph is a subset of~$\mathbb{H}\cup\mathbb{R}_\infty$, which can be viewed as a planar graph. The vertices of~$\mathcal{G}$ all belong to~$\mathbb{R}_\infty$: they are the rationals together with the point~$\infty$. From now on, we assume that every rational ~$a/b$ is in reduced form, meaning that ~$a$ and~$b$ are coprime, and~$b$ is positive. The edges of~$\mathcal{G}$ are hyperbolic geodesics in~$\mathbb{H}$: two rationals~$a/b$ and~$c/d$ are joined by an edge of~$\mathcal{G}$ if and only if~$|ad-bc|=1$ (with the convention that~$\infty$ is identified with~$1/0$).  The Farey graph induces a tessellation of the hyperbolic plane that also appears in \cite{Si1988}, as a universal cover for any triangular map on a surface. Part of the Farey graph is shown in Figure~\ref{ixe} (both grey and black lines).

\begin{figure}[ht]
\centering
\includegraphics{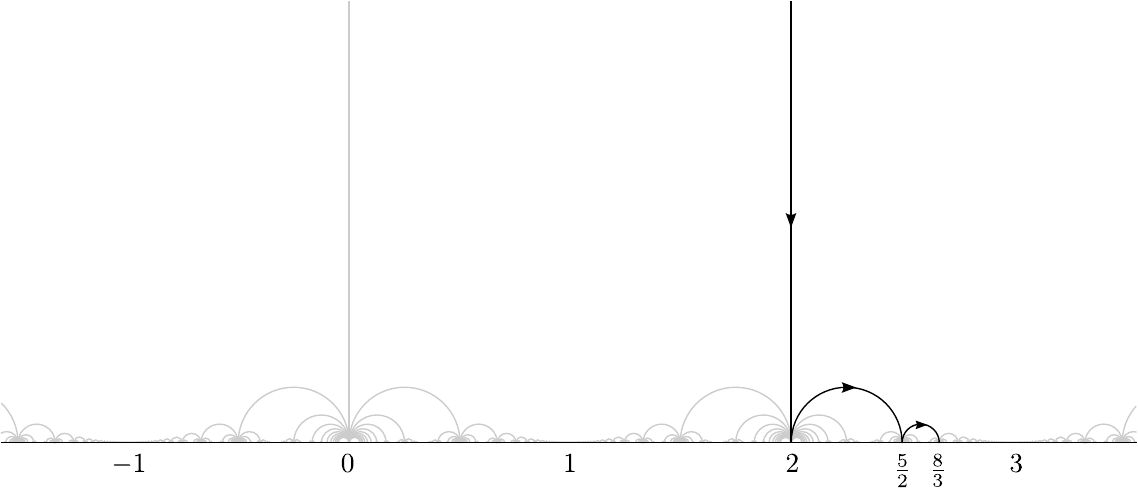}
\caption{A path in the Farey tree}
\label{figure2}
\end{figure}

The \emph{Farey tree}, which we denote by~$\mathcal{F}$, is obtained by removing all vertices  that as rationals in reduced form have odd numerator and denominator. It is a tree with a countably infinite number of vertices, and countably infinitely many edges incident to each vertex. The vertices adjacent to~$\infty$ are the even integers. Part of the Farey tree is shown in black in Figure~\ref{ixe}, and there is another illustration of~$\mathcal{F}$ in Figure~\ref{figure2} without the distraction of the Farey graph. The Farey tree induces a tessellation of the hyperbolic plane, which is Singerman's universal tessellation -- although the definition in \cite{Si1988} is slightly different to this one. (We remark that in some other works `Farey tree' refers to a different subgraph of~$\mathcal{G}$ than~$\mathcal{F}$.)  

There are other ways to define~$\mathcal{F}$ and~$\mathcal{G}$. Here is one such way. Let~$\ell$ denote the hyperbolic geodesic in~$\mathbb{H}$ between~$0$ and~$\infty$. Then the edges of~$\mathcal{G}$ are the images of~$\ell$ under the modular group~$\Gamma$ (and the vertices of~$\mathcal{G}$ are the images of~$\infty$ under~$\Gamma$). We can describe~$\mathcal{F}$ in a similar manner. Let~$\Theta$ denote the group generated by the transformations~$s(z)=-1/z$ and~$h(z)=z+2$. This Fuchsian group, called the \emph{theta group}, is a subgroup of the modular group of index 3. It consists of those M\"obius transformations~$z\mapsto (az+b)/(cz+d)$, where~$a,b,c,d\in\mathbb{Z}$ and~$ad-bc=1$, such that
\[
\begin{pmatrix}a & b \\ c& d\end{pmatrix}\equiv \begin{pmatrix}1 & 0 \\ 0& 1\end{pmatrix}\,\text{or}\,\begin{pmatrix}0 & $1$ \\ 1& 0\end{pmatrix} \pmod{2}
\]
(see \cite[Corollary~4]{Kn1970}). The edges of~$\mathcal{F}$ are the images of~$\ell$ under~$\Theta$ (and the vertices of~$\mathcal{F}$ are the images of~$\infty$ under~$\Theta$). 

We call the vertices of the Farey tree \emph{$\infty$-rationals}. They are reduced rationals whose numerator and denominator differ in parity, together with the point~$\infty$. The~$\infty$-rationals are the fixed points of one of the two conjugacy classes of parabolic elements in~$\Theta$. The vertices of~$\mathcal{G}$ that are not vertices of~$\mathcal{F}$ are called \emph{$1$-rationals} because they consist of the images of $1$ under~$\Theta$. They are the reduced rationals with odd numerator and denominator (called \emph{face-centre points} in \cite{Si1988}), and they are the fixed points of the other of the two conjugacy classes of parabolic elements in~$\Theta$. It can easily be shown that~$\Theta$ acts on~$\mathcal{F}$, and in fact each element of~$\Theta$ is a graph automorphism of~$\mathcal{F}$.

This paper is about an attractive connection between the Farey tree and \emph{even-integer continued fractions}. An even-integer continued fraction (or, more briefly, an EICF) is a sequence of even integers~$b_1,b_2,\dotsc$, which may be finite or infinite (or empty), such that all terms except possibly~$b_1$ are nonzero. We denote this continued fraction by~$[b_1,b_2,\dotsc]$ (and sometimes by~$[b_1,\dots,b_n]$ if it is finite). The number
\[
b_1 + \cfrac{1}{b_2
          + \cfrac{1}{b_3 + \cdots+\cfrac{1}{b_n}}}\, ,
\]
is called the \emph{value} of the finite EICF~$[b_1,\dots,b_n]$. The \emph{convergents} of a finite or infinite EICF~$[b_1,b_2,\dotsc]$ are the values of~$[b_1,\dots,b_n]$ for~$n=1,2,\dotsc$. If the sequence of convergents of an infinite EICF converges in~$\mathbb{R}_\infty$ to a point~$x$, then we say that the EICF \emph{converges} and has \emph{value}~$x$. Sometimes we abuse notation and use ~$[b_1,b_2,\dotsc]$ to represent its value; this is quite natural -- in fact, the distinction between continued fractions and their values is blurred in most works on continued fractions. An EICF \emph{expansion} of a real number~$x$ is an EICF with value~$x$.

In either of the graphs~$\mathcal{F}$ or~$\mathcal{G}$, we say that two vertices are \emph{adjacent} or \emph{neighbours} if they are incident to the same edge. A \emph{path} in one of these graphs is a sequence of \emph{distinct} vertices~$v_1,v_2,\dotsc$ such that~$v_i$ and~$v_{i+1}$ are adjacent for~$i=1,2,\dotsc$. The path is said to be \emph{finite} if the sequence has finite length, and otherwise it is \emph{infinite}. We say that an infinite path~$v_1,v_2,\dotsc$ \emph{converges} to a real number~$x$ if the sequence converges to~$x$ in~~$\mathbb{R}_\infty$. In these circumstances, we describe~$v_1,v_2,\dotsc$ as a path \emph{from~$v_1$ to~$x$}.

Let
\[
t_n(z)=b_n+\frac{1}{z}\quad\text{and}\quad T_n=t_1\circ t_2 \circ \dots \circ t_n,\quad n=1,2,\dotsc,
\]
where~$b_1,b_2,\dotsc$ are even integers and all except possibly~$b_1$ are nonzero. Notice that the convergents of the EICF~$[b_1,b_2,\dotsc]$ are~$T_1(\infty),T_2(\infty),\dotsc$. Now,~$0$ and~$\infty$ are adjacent in~$\mathcal{F}$, and it is easy to check that adjacency is preserved by the maps~$t_n$, so~$T_n(0)$ and~$T_n(\infty)$ are also adjacent in~$\mathcal{F}$. But
\[
T_n(0)=T_{n-1}t_n(0)=T_{n-1}(\infty),
\]
so any two consecutive vertices in the sequence~$\infty, T_1(\infty),T_2(\infty),\dotsc$ are adjacent. Furthermore, the condition~$b_n\neq 0$ for~$n\geq 2$ implies that this walk in~$\mathcal{F}$  never `backtracks': it is a path. Conversely, a short argument shows that the vertices of a path with initial vertex~$\infty$  are the convergents of a unique EICF. Thus we see that there is a \emph{correspondence between even-integer continued fractions and paths in~$\mathcal{F}$ with initial vertex~$\infty$}. Finite continued fractions correspond to finite paths, and infinite continued fractions correspond to infinite paths (and the empty continued fraction corresponds to the path consisting of the vertex~$\infty$ alone). 

For example, the EICF expansion of the rational~$8/3$ is~$[2,2,-2]$, and this continued fraction corresponds to the path in~$\mathcal{F}$ represented by the black directed edges in Figure~\ref{figure2}. The vertices of this path are, in order,~$\infty, 2, 5/2, 8/3$ and the final three of these are the convergents of the continued fraction. 

There is a similar correspondence between \emph{integer} continued fractions and paths in the Farey \emph{graph} that is well known (and the proofs of the validity of the correspondence are similar); see, for example, \cite{BeHoSh2012,Sc2011}. However, there are two reasons why the tree~$\mathcal{F}$ is better to work with than the graph~$\mathcal{G}$: (i) all infinite paths in the tree converge, and (ii)  there is an (almost) unique path from $\infty$ to each real number (in particular, as $\mathcal{F}$ is a tree there is a unique finite path from $\infty$ to each $\infty$-rational). In terms of even-integer continued fractions, these statements are (i) all infinite EICFS converge, and (ii) each real number has an (almost) unique EICF expansion. We explain the meaning of the qualification `almost' later on. Both (i) and (ii) fail for integer continued fractions, but they do hold for \emph{regular} continued fractions (the most familiar type of continued fractions, with positive integer coefficients). Here we will show that in fact much of the theory of regular continued fractions  (from, for example, \cite[Chapters~I~and~II]{Kh1997} or \cite[Chapter~X]{HaWr2008}) can be reformulated using even-integer continued fractions. To an extent, this is already known, and has been demonstrated in works such as \cite{KrLo1996,Sc1982}. The novelty of our approach is that we develop the theory of even-integer continued fractions geometrically using elementary properties of the Farey tree. 

 In Sections~2--4 we prove some of the more fundamental theorems on even-integer continued fractions using the Farey tree, covering material that is similar (although not identical) to part of \cite{KrLo1996}. Sections~5 and~6 contain results that appear to be new. To keep this account concise, we omit certain relevant topics such as the EICF expansions of quadratic irrationals and the Hurwitz constant for the theta group (see \cite{Sc1940} for a treatment of the latter topic in the spirit of this paper). Furthermore, for the sake of brevity, we sometimes skip the details of  elementary geometric arguments, so that the reader gets a feel for the geometric approach without getting bogged down in details.

\section{Infinite continued fractions}

In this section we prove that every infinite EICF converges. There are several ways to do this; for example, we could invoke a more general theorem on the convergence of continued fractions, or we could use algebraic relationships between the convergents to estimate the distance between consecutive convergents. Our approach is to use the Farey tree to establish the following theorem. 

\begin{theorem}\label{dnm}
Every infinite EICF converges to an irrational or a~\text{$1$-rational}.
\end{theorem}

To prove the theorem, consider any infinite EICF, and let $\gamma$ be the corresponding infinite path in $\mathcal{F}$ with initial vertex $\infty$. First we will show that~$\gamma$ cannot accumulate at an $\infty$-rational. Suppose, on the contrary, that~$\gamma$ does accumulate at a vertex~$x$ of~$\mathcal{F}$. By applying an element of~$\Theta$ to~$\gamma$ if necessary we can assume that~$x\neq\infty$. Furthermore, by removing the first so many terms from~$\gamma$ we can assume that it does not pass through~$x$. Let~$a$ be the initial vertex of~$\gamma$.

\begin{figure}[ht]
\centering
\begin{tikzpicture}[scale=2.7]
\draw (0.5,0) arc (0:180:0.5);
\draw (1,0) arc (0:180:0.25);
\draw  (-1.5,0) -- (1.5,0);
\filldraw (-1,0) circle (0.6pt) node[below]{$a$};
\filldraw (-0.5,0) circle (0.6pt) node[below]{$u$};
\filldraw (0.5,0) circle (0.6pt) node[below]{$x$};
\filldraw (1,0) circle (0.6pt) node[below]{$v$};
\end{tikzpicture}
\caption{Two neighbours~$u$ and~$v$ of the vertex~$x$, and another vertex~$a$}
\label{wuj}
\end{figure}
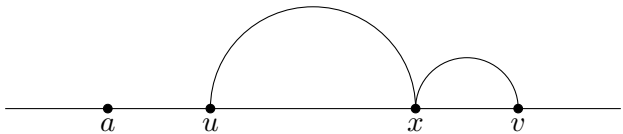

Like all vertices of~$\mathcal{F}$, the vertex~$x$ has infinitely many neighbours, which accumulate on the left and right of~$x$. Choose any two neighbours~$u$ and~$v$ such that~$u<x<v$, and such that the vertex~$a$ lies outside the real interval~$(u,v)$, as shown in Figure~\ref{wuj}. Edges of~$\mathcal{F}$ do not intersect in~$\mathbb{H}$, so we see that because~$\gamma$ accumulates at~$x$, it must pass through one of~$u$,~$x$ and~$v$. However, because~$\mathcal{F}$ is a tree, with only one exception, any path from~$a$ to a neighbour of~$x$ must pass through~$x$. So providing we choose~$u$ and~$v$ such that neither of them are the exceptional neighbour, we see that~$\gamma$ passes through~$x$. This contradicts an earlier assumption, so ~$\gamma$ cannot accumulate at a vertex of~$\mathcal{F}$ after all. 

We have just seen that the path $\gamma$ cannot accumulate at an $\infty$-rational. Suppose, in order to reach a contradiction, that~$\gamma$ accumulates at two numbers~$x$ and~$y$, each of which is either irrational or a~$1$-rational, and~$x<y$. Now, the vertices of~$\mathcal{F}$ that lie inside the real interval~$(x,y)$ are connected in~$\mathcal{F}$ to the vertices that lie outside this interval, so in particular there must be an edge  of~$\mathcal{F}$ with one end vertex $u$ inside the interval and the other $v$ outside. Edges of~$\mathcal{F}$ do not intersect in~$\mathbb{H}$, so we see that because~$\gamma$ accumulates at both~$x$ and~$y$, it must pass through at least one of $u$ or $v$ infinitely many times, which is impossible. Thus, contrary to our assumption,~$\gamma$ cannot accumulate at two numbers, so it converges. The proof of Theorem~\ref{dnm} is now complete.

\section{Representing real numbers by even-integer continued fractions}

The next fundamental result is about the existence and uniqueness of EICF expansions of real numbers.  It is unoriginal (see, for example, \cite{KrLo1996}, where there are a number of results similar to parts of this one); however, our method of proof using the Farey tree is original, and it is simple and elegant. \newpage

\begin{theorem}\label{ixy}\quad
\begin{enumerate}[label=\emph{(\roman*)},leftmargin=20pt,topsep=0pt]
\item The value of any finite EICF is an~$\infty$-rational, and each~$\infty$-rational has a unique finite EICF expansion.
\item The value of an infinite EICF is either  irrational or a~$1$-rational, and
\begin{enumerate}[label=\emph{(\alph*)}]
\item each irrational has a unique infinite EICF expansion, 
\item each $1$-rational has exactly two infinite EICF~expansions, each of which eventually alternates between~$2$ and~$-2$. 
\end{enumerate}
\end{enumerate}
\end{theorem}

 As~$\mathcal{F}$ is a tree, and the vertices are the~$\infty$-rationals, we can immediately deduce statement (i) of the theorem using the correspondence between even-integer continued fractions and paths in~$\mathcal{F}$. We now turn to statement (ii). The first part of statement (ii) follows from Theorem~\ref{dnm}. It remains only to discuss statements (a) and (b). 

\begin{figure}[ht]
\centering
\includegraphics{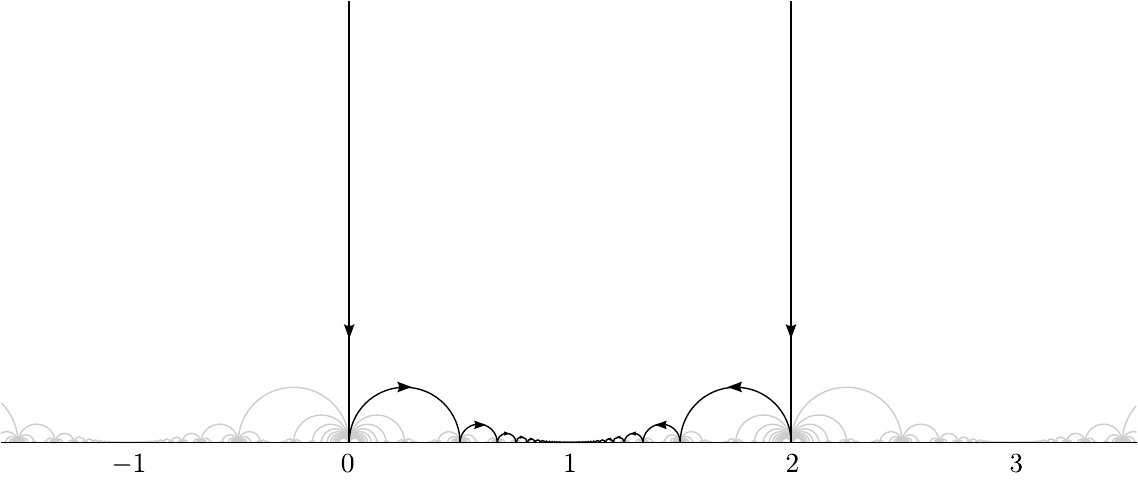}
\caption{Two paths that converge to 1}
\label{figure4}
\end{figure}

We begin this discussion by looking at EICF expansions of the number 1; here are two of them:
\[
1=[0,2,-2,2,-2,\dotsc]=[2,-2,2,-2,\dotsc].
\]
We can check that the value~$x$ of the second continued fraction is $1$ by observing that~$x$ must satisfy
\[
x = 2+\cfrac{1}{-2+\cfrac{1}{x}},
\]
and the only solution of this equation is~$x=1$. (The value of the first continued fraction can be obtained in a similar manner.) The paths in~$\mathcal{F}$ corresponding to these two continued fractions are shown marked by arrows in Figure~\ref{figure4}.

\begin{figure}[ht]
\centering
\includegraphics{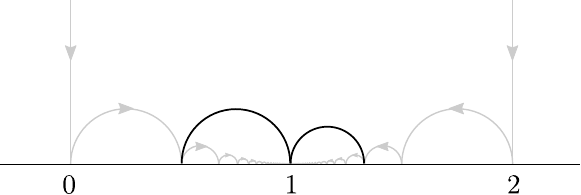}
\caption{Two neighbours of 1}
\label{figure5}
\end{figure}

In fact, the two EICF expansions that we have found are the \emph{only} EICF expansions of 1. To see why this is so, observe that in the Farey graph~$\mathcal{G}$, every single one of the vertices in these two paths is connected to $1$ by an edge (in fact, they are the full collection of neighbours of $1$ in~$\mathcal{G}$ -- see Figure~\ref{ixe}). Two such edges are shown in Figure~\ref{figure5}, on either side of 1.

Suppose now that~$\gamma$ is an infinite path in~$\mathcal{F}$ from~$\infty$ to 1. Aside from the initial vertex~$\infty$, this path must lie entirely to the left or entirely to the right of $1$ (because any path in~$\mathcal{F}$ that passes from one side to the other of~$1$ must pass through~$\infty$). Suppose that it lies to the left -- the other case can be handled in a similar way. Then because edges in the Farey graph do not intersect,~$\gamma$ must pass through all of the vertices of~$\alpha$. There is only one such path that does this, namely~$\alpha$ itself, so~$\gamma=\alpha$.

We summarise this discussion in a lemma.

\begin{lemma}\label{yyv}
The number $1$ has precisely two EICF expansions, namely $[0,2,-2,2,-2,\dotsc]$ and~$[2,-2,2,-2,\dotsc]$.
\end{lemma}

If~$x$ is any $1$-rational, then there is an element~$g$ of~$\Theta$ such that~$g(1)=x$. It follows that~$g(\alpha)$ and~$g(\beta)$ are infinite  paths from~$g(\infty)$ to 1. By connecting~$\infty$ to~$g(\infty)$ we obtain two walks from~$\infty$ to $1$ (each  may have repeated vertices), which we can modify by adjusting a finite number of terms to give two paths from $\infty$ to 1. Thus we obtain two EICF expansions of~$x$. We can reverse this argument to see that these are the only EICF expansions of~$x$. This gives us the following corollary of Lemma~\ref{yyv}.

\begin{corollary}\label{qlz}
Every $1$-rational has precisely two EICF expansions.
\end{corollary}

In the next section we will see that if $x$ and $y$ have infinite EICF expansions, and $g(x)=y$ for some transformation $g$ in $\Theta$, then it is possible to remove a finite number of consecutive terms from the start of the EICF expansions of $x$ and $y$ to give two expansions that agree. Therefore the $1$-rationals are precisely those real numbers that have an EICF expansion that eventually alternates between $2$ and $-2$. Furthermore, it is straightforward to check that the two continued fractions
\[
[b_1,\dots,b_n,2,-2,2,\dotsc]\quad\text{and}\quad [b_1,\dots,b_{n-1},b_n+2,-2,2,-2,\dotsc]
\]
have the same value, so the two EICF expansions referred to in Corollary~\ref{qlz} are of these forms.

We have now proved statement (b) of Theorem~\ref{ixy}, which leaves only statement (a). Let us prove the uniqueness assertion of (a). Suppose then that~$\alpha$ and~$\beta$ are two infinite paths from~$\infty$ to a real number~$x$. The two paths may coincide for a certain number of vertices: let~$w$ be the final vertex for which they do so. Choose an element~$g$ of~$\Theta$ such that~$g(w)=\infty$. Let~$\alpha'$ and~$\beta'$ be the paths obtained from~$g(\alpha)$ and~$g(\beta)$, respectively, after removing all vertices that occur before~$\infty$. Then~$\alpha'$ and~$\beta'$ are infinite paths from~$\infty$ to~$g(x)$, such that the second vertex~$u$ of~$\alpha'$ is distinct from the second vertex~$v$ of~$\beta'$. The vertices~$u$ and~$v$ are even integers, so there is an odd integer~$q$ (a $1$-rational) that lies between them on the real line. Neither~$\alpha'$ nor~$\beta'$ can cross~$q$, and since they converge to the same value, that value must be~$q$.  Therefore~$g(x)$ is a $1$-rational, so~$x$ is also a $1$-rational. 

This argument shows that each irrational has at most one EICF expansion. Let us now show that each irrational has at least one such expansion. One way to do this is to use an algorithm of a similar type to Euclid's algorithm: in this case the `nearest even-integer algorithm' does the trick. However, we prefer to justify the existence of an expansion using the Farey graph and tree. 

We define a \emph{Farey interval} to be a real interval whose endpoints are neighbouring vertices in the Farey graph~$\mathcal{G}$. If~$[a/b,c/d]$ is a Farey interval (where, as usual, the fractions are given in reduced form), then it is easily seen that~$[a/b,(a+c)/(b+d)]$ and~$[(a+c)/(b+d),b/d]$ are both Farey intervals -- let us call them the \emph{Farey subintervals} of~$[a/b,c/d]$. Now, any irrational~$x$ belongs to a Farey interval~$[n,n+1]$, where~$n$ is the integer part of~$x$, and by repeatedly choosing Farey subintervals, we can construct a nested sequence of Farey intervals that contains~$x$ in its intersection.  The width of one of these intervals~$[a/b,c/d]$ is 
\[
\left|\frac{a}{b}-\frac{c}{d}\right| = \left|\frac{ad-bc}{bd}\right|=\frac{1}{bd},
\]
so we see that the sequence of widths of this nested sequence of Farey intervals converges to 0. 

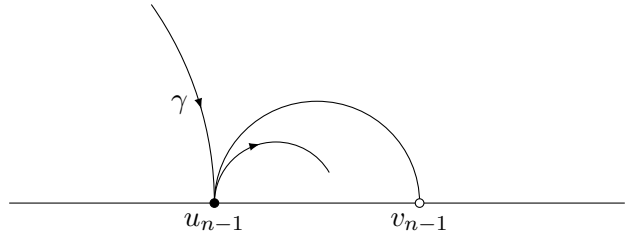
\begin{figure}[ht]
\centering
\begin{tikzpicture}[scale=2.7]
\draw (0.5,0) arc (0:180:0.5);
\draw  (-1.5,0) -- (1.5,0);
\filldraw (-0.5,0) circle (0.6pt) node[below]{$u_{n-1}$};
\draw[fill=white] (0.5,0) circle (0.6pt) node[below]{$v_{n-1}$};
\draw[middlearrow={latex reversed}] (-0.5,0) node[above, yshift=30pt, xshift=-13pt] {$\gamma$} arc (0:35:1.7);
\draw[middlearrow={latex}] (-0.5,0) arc (180:30:0.3);
\end{tikzpicture}
\caption{The path $\gamma$ passes through $u_{n-1}$}
\label{poa}
\end{figure}

Let us now restrict attention to those infinitely many Farey intervals $I_1\supset I_2\supset\dotsb$ from the sequence for which one of the endpoints of~$I_n$ is a $1$-rational~$v_n$ (and the other endpoint~$u_n$ must then be an~$\infty$-rational). Let~$\gamma_n$ be the unique path from~$\infty$ to~$u_n$ in~$\mathcal{F}$. Any path in~$\mathcal{F}$ from~$\infty$ to a vertex inside $I_{n-1}$ must pass through~$u_{n-1}$ (because $u_{n-1}$ and $v_{n-1}$ are neighbours in $\mathcal{G}$, as illustrated in Figure~\ref{poa}, and edges of $\mathcal{G}$ do not intersect). Therefore~$\gamma_{n-1}$ is a subpath of~$\gamma_n$. It follows that there is a unique infinite path~$\gamma$ that contains every path~$\gamma_n$ as a subpath. The path~$\gamma$ passes through all the vertices~$u_n$, which accumulate at~$x$, so~$\gamma$ must converge to~$x$. This completes the proof of Theorem~\ref{ixy}.

\section{Serret's theorem on continued fractions}

This section is about a counterpart for even-integer continued fractions of a well-known theorem of Serret on regular continued fractions. Before we state our theorem, we must introduce the \emph{extended theta group}, which is the group~$\widetilde\Theta$ generated by the theta group and the transformation~$r(z)=-z$. This group acts on~$\mathbb{R}_\infty$, and it also acts on the set of~$\infty$-rationals. In fact,  elements of~$\widetilde\Theta$ preserve adjacency in~$\mathcal{F}$, so~$\widetilde\Theta$ acts on the abstract graph underlying~$\mathcal{F}$. We say that two real numbers are \emph{equivalent} under the action of $\widetilde\Theta$ if they lie in the same orbit under this action.

Our version of Serret's theorem for even-integer continued fractions follows. It is similar to \cite[Theorem~1]{KrLo1996}, but not quite the same  because even-integer continued fractions are defined differently in that paper.

\begin{theorem}\label{dwl}
Two real numbers~$x$ and~$y$ that are not~$\infty$-rationals are equivalent under~$\widetilde\Theta$ if and only if there are positive integers~$m$ and~$n$ such that the EICF expansions of~$x$ and~$y$,
\[
x=[a_1,a_2,\dotsc]\quad\text{and}\quad y=[b_1,b_2,\dotsc],
\]
either satisfy~$a_{m+i}=b_{n+i}$ for~$i=1,2,\dotsc$ or~$a_{m+i}=-b_{n+i}$ for~$i=1,2,\dotsc$. 
\end{theorem}

Serret's theorem for regular continued fraction expansions is similar, but uses an extension of the modular group rather than the theta group, and the condition~$a_{m+i}=-b_{n+i}$ for~$i=1,2,\dotsc$ is absent.

Crucial to the proof of this theorem is the following lemma.

\begin{lemma}\label{tab}
If a real number~$x$ has an EICF expansion~$[a_1,a_2,\dotsc]$, then an EICF expansion of~$-x$ is~$[-a_1,-a_2,\dotsc]$.
\end{lemma}

There is no obvious analogue of this lemma for regular continued fractions because the coefficients of regular continued fractions are (almost) all positive. 

The lemma can be proven with the Farey tree by observing that the paths from~$\infty$ to~$x$ and from~$\infty$ to~$-x$ are reflections of each other in the imaginary axis. However, in this case, will prove the lemma using M\"obius transformations. Let~$t_a(z)=a+1/z$, where~$a$ is even; this transformation belongs to~$\widetilde\Theta$. Observe that~$rt_ar=t_{-a}$. We are given that an EICF expansion of~$x$ is~$[a_1,a_2,\dotsc]$, which implies that~$t_{a_1}t_{a_2}\dotsb t_{a_n}(\infty)\to~x$ as~$n\to\infty$. Now
\[
t_{-a_1}t_{-a_2}\dotsb t_{-a_n}(\infty)=rt_{a_1}t_{a_2}\dotsb t_{a_n}r(\infty)=rt_{a_1}t_{a_2}\dotsb t_{a_n}(\infty).
\]
So~$t_{-a_1}t_{-a_2}\dotsb t_{-a_n}(\infty)\to r(x)=-x$ as~$n\to\infty$. Therefore an EICF expansion of~$-x$ is~$[-a_1,-a_2,\dotsc]$.

Let us now prove Theorem~\ref{dwl}. Suppose first that~$y=g(x)$, where~$g\in\widetilde\Theta$. We wish to prove that there are positive integers~$m$ and~$n$ such that~$a_{m+i}=b_{n+i}$ for~$i=1,2,\dotsc$ or~$a_{m+i}=-b_{n+i}$ for~$i=1,2,\dotsc$. Since~$\widetilde\Theta$ is generated by the transformations~$r(z)=-z$,~$t(z)=1/z$ and~$h(z)=z+2$, it suffices to prove the assertion when~$g$ is each of~$r$,~$t$,~$h$ and~$h^{-1}$. It is straightforward to do so when~$g$ is one of the final three transformations, and the remaining case when~$g$ is~$r$ is an immediate consequence of Lemma~\ref{tab}. 

For the converse, suppose  that~$x=[a_1,a_2,\dotsc]$,~$y=[b_1,b_2,\dotsc]$ and either (i)~$a_{m+i}=b_{n+i}$ for~$i=1,2,\dotsc$, or (ii)~$a_{m+i}=-b_{n+i}$ for~$i=1,2,\dotsc$. By replacing~$x$ by~$-x$ if necessary, and invoking Lemma~\ref{tab}, we can assume that (i) holds. Observe that
\[
x = t_{a_1}\dotsb t_{a_m}([a_{m+1},a_{m+2},\dotsc])\quad \text{and}\quad y = t_{b_1}\dotsb t_{b_n}([b_{n+1},b_{n+2},\dotsc]).
\]
Hence~$y=t_{b_1}\dotsb t_{b_n}t_{a_m}^{-1}\dotsb t_{a_1}^{-1}(x)$, so~$x$ and~$y$ are equivalent under~$\widetilde\Theta$. This completes the proof of Theorem~\ref{dwl}.

\section{An alternative characterisation of convergents of even-integer continued fractions}

In this section we describe an alternative way to characterise the convergents of the EICF expansion of any irrational~$x$. The characterisation can easily be adapted to allow~$x$ to be rational. 

\begin{theorem}\label{icd}
A finite $\infty$-rational~$u$  is a convergent of the EICF expansion of an irrational~$x$ if and only if there is a $1$-rational~$v$ adjacent to~$u$ in the Farey graph such that~$x$ lies between~$u$ and~$v$ on the real line.
\end{theorem}

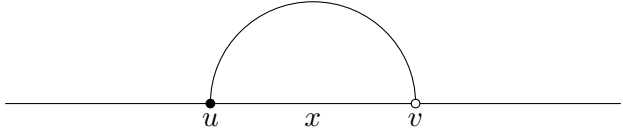
\begin{figure}[ht]
\centering
\begin{tikzpicture}[scale=2.7]
\draw (0.5,0) arc (0:180:0.5);
\draw  (-1.5,0) -- (1.5,0);
\filldraw (-0.5,0) circle (0.6pt) node[below]{$u$};
\draw[fill=white] (0.5,0) circle (0.6pt) node[below]{$v$};
\node[below] at (0,0){$x$};
\end{tikzpicture}
\caption{The irrational $x$ lies between the $\infty$-rational $u$ and the $1$-rational $v$}
\label{cie}
\end{figure}

The second statement of the theorem is illustrated in Figure~\ref{cie}. Its significance for the Farey graph~$\mathcal{G}$ is that any vertex~$w$ of~$\mathcal{G}$ that is sufficiently close to~$x$ on the real line is separated from~$\infty$ by the edge incident to~$u$ and~$v$. So any path from~$\infty$ to~$w$ must pass through one of~$u$ or~$v$ -- and if the path lies in~$\mathcal{F}$, then it must pass through~$u$. In particular, this demonstrates that~$u$ must be a convergent of the EICF expansion of~$x$.

The converse implication of Theorem~\ref{icd} is a direct consequence of the following lemma (which is a slightly stronger statement).

\begin{lemma}\label{lak}
Let~$u$ and~$w$ be two consecutive convergents in the EICF expansion of an irrational~$x$, in that order. Then there is a $1$-rational~$v$ adjacent to each of~$u$ and~$w$ in the Farey graph such that both~$w$ and~$x$ lie between~$u$ and~$v$ on the real line. 
\end{lemma}

\begin{figure}[ht]
\centering
\begin{tikzpicture}[scale=2.7]
\draw  (-1.5,0) -- (1.5,0);
\draw (0.5,0) arc (0:180:0.5);
\draw (0,0) arc (0:180:0.25);
\draw (0.5,0) arc (0:180:0.25);
\filldraw (-0.5,0) circle (0.6pt) node[below]{$u$};
\filldraw (-0,0) circle (0.6pt) node[below]{$w$};
\draw[fill=white] (0.5,0) circle (0.6pt) node[below]{$v$};
\end{tikzpicture}
\caption{A triangle in the Farey graph}
\label{iio}
\end{figure}
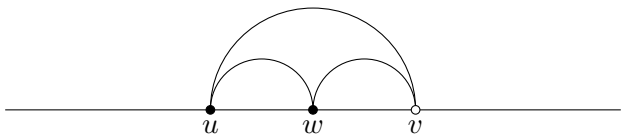

Since~$u$ and~$w$ are adjacent in~$\mathcal{F}$, they are also adjacent in~$\mathcal{G}$. There are two other vertices in~$\mathcal{G}$ that are adjacent to both~$u$ and~$w$, precisely one of which (call it~$v$) does \emph{not} lie between~$u$ and~$w$ on the real line.  Let $\gamma$ be the path of convergents of the EICF expansion of~$x$. If $\gamma$ enters the interval between $u$ and $v$, then it must pass through $u$ to get there, and it cannot leave the interval. Similar comments apply to the interval between~$w$ and~$v$. Now,  $u$ cannot lie in the interval between $w$ and $v$  because if it did then, as we have just seen, the path $\gamma$ would pass through $w$ before it passed through~$u$. So~$w$ lies in the  interval between $u$ and $v$ (as illustrated in Figure~\ref{iio}), and~$x$ lies in that interval too. This completes the proofs of Lemma~\ref{lak} and Theorem~\ref{icd}.

\section{Approximating irrationals by rationals}

One of the principal uses of continued fractions is in the field of Diophantine approximation, which is concerned with approximating real numbers by rationals. In this section we prove an analogue for even-integer continued fractions of a classic result of Lagrange on regular continued fractions.

We call an~$\infty$-rational~$a/b$ a \emph{strong~$\infty$-approximant} of a real number~$x$ if for each~$\infty$-rational~$c/d$ such that~$d\leq b$, we have
\[
|bx-a|\leq|dx-c|,
\]
with equality if and only if~$c/d=a/b$. 

\begin{theorem}\label{ood}
An~$\infty$-rational is a strong~$\infty$-approximant of an irrational~$x$ if and only if it is a convergent of the EICF expansion of~$x$.
\end{theorem}

Lagrange's theorem for regular continued fractions is similar (see \cite[Theorems~16~and~17]{Kh1997}), but uses rationals rather than~$\infty$-rationals. 

There is no need to assume that~$x$ is irrational in the theorem -- subject to minor modifications of the theorem we can allow $x$ to be any real number -- but it is the irrational case that interests us most, and the proof is marginally simpler with the assumption that~$x$ is irrational.

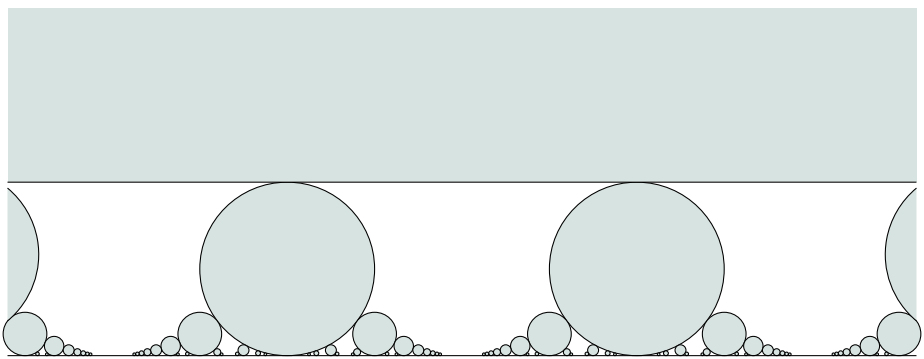
\begin{figure}[ht]
\centering
\begin{tikzpicture}[scale=2.3]
\horo{0}{1/2}{col1}\horo{2}{1/2}{col1}

\horo{-3/2}{1/8}{col1}\horo{-1/2}{1/8}{col1}\horo{1/2}{1/8}{col1}\horo{3/2}{1/8}{col1}\horo{5/2}{1/8}{col1}\horo{7/2}{1/8}{col1}

\horo{-4/3}{1/18}{col1}\horo{-2/3}{1/18}{col1}\horo{2/3}{1/18}{col1}\horo{4/3}{1/18}{col1}\horo{8/3}{1/18}{col1}\horo{10/3}{1/18}{col1}
  
\horo{-5/4}{1/32}{col1}\horo{-3/4}{1/32}{col1}\horo{-1/4}{1/32}{col1}\horo{1/4}{1/32}{col1}\horo{3/4}{1/32}{col1}\horo{5/4}{1/32}{col1}\horo{7/4}{1/32}{col1}\horo{9/4}{1/32}{col1}\horo{11/4}{1/32}{col1}\horo{13/4}{1/32}{col1}

\horo{-6/5}{1/50}{col1}\horo{-4/5}{1/50}{col1}\horo{-2/5}{1/50}{col1}\horo{2/5}{1/50}{col1}\horo{4/5}{1/50}{col1}\horo{6/5}{1/50}{col1} \horo{8/5}{1/50}{col1}\horo{12/5}{1/50}{col1}\horo{14/5}{1/50}{col1}\horo{16/5}{1/50}{col1}

\horo{-7/6}{1/72}{col1}\horo{-5/6}{1/72}{col1}\horo{-1/6}{1/72}{col1}\horo{1/6}{1/72}{col1}\horo{5/6}{1/72}{col1}\horo{7/6}{1/72}{col1}\horo{11/6}{1/72}{col1}\horo{13/6}{1/72}{col1}\horo{17/6}{1/72}{col1}\horo{19/6}{1/72}{col1}

\horo{-10/7}{1/98}{col1}\horo{-8/7}{1/98}{col1}\horo{-6/7}{1/98}{col1}\horo{-4/7}{1/98}{col1}\horo{-2/7}{1/98}{col1}\horo{2/7}{1/98}{col1}\horo{4/7}{1/98}{col1}\horo{6/7}{1/98}{col1}\horo{8/7}{1/98}{col1}\horo{10/7}{1/98}{col1}\horo{12/7}{1/98}{col1}\horo{16/7}{1/98}{col1}\horo{18/7}{1/98}{col1}\horo{20/7}{1/98}{col1}\horo{22/7}{1/98}{col1}\horo{24/7}{1/98}{col1}

\horo{-11/8}{1/128}{col1}\horo{-9/8}{1/128}{col1}\horo{-7/8}{1/128}{col1}\horo{-5/8}{1/128}{col1}\horo{-3/8}{1/128}{col1}\horo{-1/8}{1/128}{col1}\horo{1/8}{1/128}{col1}\horo{3/8}{1/128}{col1}\horo{5/8}{1/128}{col1}\horo{7/8}{1/128}{col1}\horo{9/8}{1/128}{col1}\horo{11/8}{1/128}{col1}\horo{13/8}{1/128}{col1}\horo{15/8}{1/128}{col1}\horo{17/8}{1/128}{col1}\horo{19/8}{1/128}{col1}\horo{21/8}{1/128}{col1}\horo{23/8}{1/128}{col1}\horo{25/8}{1/128}{col1}\horo{27/8}{1/128}{col1}

\draw [fill=col1,draw=none] (-1.6,1) rectangle (3.6,2);
\draw [fill=col1,draw=none] (-1.6,0.2) -- (-1.6,0.2) arc(-50:50:1/2) -- cycle;	
\draw (-1.6,0.2) arc(-50:50:1/2);
\draw [fill=col1,draw=none] (3.6,0.2) -- (3.6,0.2) arc(230:130:1/2) -- cycle;	
\draw (3.6,0.2) arc(230:130:1/2);
\draw (-1.6,1) -- (3.6,1);
\draw  (-1.6,0) -- (3.6,0);
\end{tikzpicture}
\caption{Ford circles based at the~$\infty$-rationals}
\label{stb}
\end{figure}

Our proof uses \emph{Ford circles}, and is similar to the proof of Lagrange's theorem from \cite{Sh2011}. Ford circles are a collection of horocycles in~$\mathbb{H}$ used by Ford to study continued fractions in papers such as \cite{Fo1917,Fo1938}. We say that a horocycle is \emph{based} at an element~$x$ of~$\mathbb{R}_\infty$ if the horocycle is tangent to~$\mathbb{R}_\infty$ at~$x$. Given a reduced rational~$u=a/b$, the Ford circle~$C_u$ is the horocycle based at~$u$ with Euclidean radius~$\text{rad}[C_u]=1/(2b^2)$. There is one other Ford circle~$C_\infty$, which is the line~$y=1$ together with the point~$\infty$. Two Ford circles intersect in at most a single point, and the interiors of the two circles are disjoint. In fact, one can check that the Ford circles~$C_{a/b}$ and~$C_{c/d}$ are tangent if and only if~$|ad-bc|=1$. Therefore the full collection of Ford circles is a model of the abstract graph underlying the Farey graph: the vertices of this graph are represented by Ford circles, and two vertices are adjacent if and only if the Ford circles are tangent. Similarly, the collection of Ford circles based at~$\infty$-rationals is a model of the abstract graph underlying the Farey tree; this model is illustrated in Figure~\ref{stb}. When studying even-integer continued fractions, it is helpful to consider both the Farey tree and this alternative model of the tree using Ford circles.

We now relate Ford circles to strong~$\infty$-approximants. Let~$u=a/b$. Notice that if~$v=c/d$, then~$d\leq b$ if and only if~$\text{rad}[C_u]\leq\text{rad}[C_v]$. For any real number~$x$, let
\[ 
R_u(x)=\frac{1}{2} |bx-a|^2.
\]
Using elementary geometry, it can be shown that~$R_u(x)$ is the Euclidean radius of the horocycle based at~$x$ that is externally tangent to~$C_u$. With this terminology, we can describe a strong~$\infty$-approximant of a real number~$x$ as an~$\infty$-rational~$u$ such that for each~$\infty$-rational~$w$ with~$\text{rad}[C_u]\leq\text{rad}[C_w]$, we have~$R_u(x)\leq R_w(x)$, with equality if and only if~$w=u$. We will use this definition of strong~$\infty$-approximants together with Theorem~\ref{icd} to prove Theorem~\ref{ood}. Our proof omits several elementary geometric details.

\begin{figure}[ht]
\centering
\begin{tikzpicture}[scale=2.7]
\draw  (-1.5,0) -- (1.5,0);
\horo{-1}{0.45}{col1}; \node[below] at (-1,0) {$w$};
\horo{-0.1}{0.3}{col1}; \node[below] at (-0.1,0) {$u$}; 
\horo{0.7}{0.5333}{col3}; \node[below] at (0.7,0) {$v$};
\draw[dashed] (0.2,0.8) circle (0.8); \node[below] at (0.2,0) {$x$};
\draw[dashed] (0.2,0.0725) circle (0.0725);
\end{tikzpicture}
\caption{Ford circles based at~$u$,~$v$ and~$w$ and horocycles based at~$x$}
\label{ujj}
\end{figure}
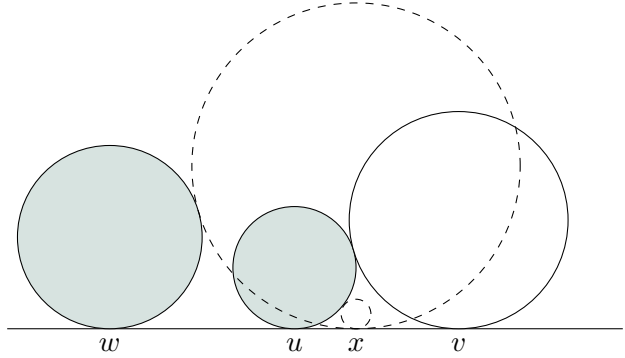

Suppose first that~$u$ is a convergent of the EICF expansion of~$x$. Theorem~\ref{icd} tells us that there is a $1$-rational~$v$ adjacent to~$u$ in the Farey graph such that~$x$ lies between~$u$ and~$v$ on the real line. If~$w$ is an~$\infty$-rational distinct from~$u$ with~$\text{rad}[C_u]\leq\text{rad}[C_w]$, then~$w$ must lie outside the real interval between~$u$ and~$v$, so~$R_u(x)< R_w(x)$, as illustrated in Figure~\ref{ujj}. Therefore~$u$ is a strong~$\infty$-approximant of~$x$.

\begin{figure}[ht]
\centering
\begin{tikzpicture}[scale=2.7]
\draw  (-1.5,0) -- (1.5,0);
\horo{-1}{5/12}{col1}; \node[below] at (-1,0) {$w_n$};
\horo{-6/11}{15/121}{col1}; \node[below] at (-6/11,0) {$w_{n+1}$}; 
\horo{0}{0.6}{col3}; \node[below] at (0,0) {$v$};
\horo{1}{0.3}{col1}; \node[below] at (1,0) {$u$};
\node[below] at (-0.23,0) {$x$};
\end{tikzpicture}
\caption{Ford circles based at~$u$,~$v$,~$w_n$ and~$w_{n+1}$}
\label{qoq}
\end{figure}
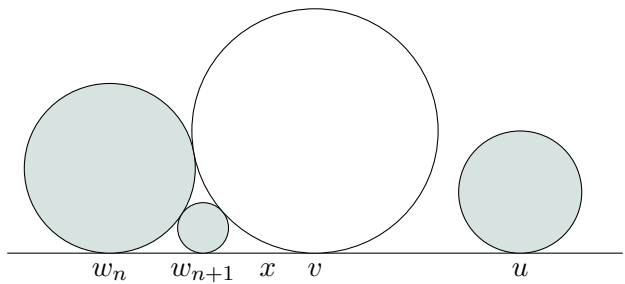

Conversely, suppose that~$u$ is an~$\infty$-rational that is not one of the convergents~$w_1,w_2,\dotsc$ of the EICF expansion of~$x$. Choose a convergent~$w_n$ such that ~$\text{rad}[C_{w_{n+1}}]< \text{rad}[C_u]\leq\text{rad}[C_{w_{n}}]$. By Theorem~\ref{icd}, there is a $1$-rational~$v$ adjacent to each of~$w_n$ and~$w_{n+1}$ in the Farey graph such that both~$w_{n+1}$ and~$x$ lie between~$w_n$ and~$v$ on the real line. On the other hand, the radius of~$C_u$ is larger than that of~$C_{w_{n+1}}$, so~$u$ does not lie between~$w_n$ and~$v$, as illustrated in Figure~\ref{qoq}. Therefore~$R_{w_n}(x)<R_u(x)$, so~$u$ is not a strong~$\infty$-approximant of~$x$. This completes the proof of Theorem~\ref{ood}.

\section{Concluding remark}

We have seen that a good deal of the theory of even-integer continued fractions can be understood by viewing such continued fractions as paths in the Farey tree. It may be of interest to study paths in other maps on surfaces, and investigate their relationship with  continued fractions.


\begin{thebibliography}{AW1}

\providecommand{\href}[2]{#2}

\bibitem{BeHoSh2012} A. F. Beardon, M. Hockman\ and\ I. Short, Geodesic continued fractions, Michigan Math. J. {\bf 61} (2012), no.~1, 133--150.

\bibitem{Fo1917} L. R. Ford, A geometrical proof of a theorem of Hurwitz, Proc. Edinburgh Math. Soc. {\bf 35} (1917), 59--65.

\bibitem{Fo1938} L. R. Ford, Fractions, Amer. Math. Monthly {\bf 45} (1938), no.~9, 586--601.

\bibitem{HaWr2008} G. H. Hardy\ and\ E. M. Wright, {\it An introduction to the theory of numbers}, sixth edition, Oxford Univ. Press, Oxford, 2008. 

\bibitem{Kh1997} A. Ya. Khinchin, {\it Continued fractions}, translated from the third (1961) Russian edition, reprint of the 1964 translation, Dover, Mineola, NY, 1997.

\bibitem{Kn1970} M. I. Knopp, {\it Modular functions in analytic number theory}, Markham Publishing Co., Chicago, IL, 1970.

\bibitem{KrLo1996} C. Kraaikamp\ and\ A. Lopes, The theta group and the continued fraction expansion with even partial quotients, Geom. Dedicata {\bf 59} (1996), no.~3, 293--333. 

\bibitem{Sc2011} R. E. Schwartz, {\it Mostly surfaces}, Student Mathematical Library, 60, Amer. Math. Soc., Providence, RI, 2011.

\bibitem{Sc1982} F.  Schweiger, Continued fractions with odd and even partial quotients, Arbeitsber. Math. Inst. Univ. Salzburg {\bf 4} (1982), 59--70.


\bibitem{Sc1940} W. T. Scott, Approximation to real irrationals by certain classes of rational fractions, Bull. Amer. Math. Soc. {\bf 46} (1940), 124--129. 

\bibitem{Sh2011} I. Short, Ford circles, continued fractions, and rational approximation, Amer. Math. Monthly {\bf 118} (2011), no.~2, 130--135.

\bibitem{Si1988} D. Singerman, Universal tessellations, Rev. Mat. Univ. Complut. Madrid {\bf 1} (1988), no.~1-3, 111--123.

\end{thebibliography}
\end{document}